\documentclass[10pt,journal,onecolumn,web]{ieeecolor}

\usepackage{generic}
\usepackage{graphicx}      
\usepackage{tabulary}
\usepackage{amsfonts}
\usepackage{amsmath}
\usepackage{amssymb}
\usepackage{pgf,tikz}
\usepackage{tabularx}
\usepackage{float}
\usepackage{cite}
\usepackage{stfloats}

\let\labelindent\relax
\usepackage{enumitem} 
\usepackage{bbm} 
\usepackage{bm}

\newtheorem{remark}{Remark}
\newtheorem{theorem}{Theorem}
\newtheorem{lemma}{Lemma}
\newtheorem{assumption}{Assumption}

\newtheorem{corollary}{Corollary}

\newcommand{\col}{{\rm col\;}}
\newcommand{\diag}{{\rm diag\;}}
\newcommand{\rank}{{\rm rank\ } }
\newcommand{\image}{{\rm image\ } }

\def\qed{ \rule{.1in}{.1in}}

\def\BibTeX{{\rm B\kern-.05em{\sc i\kern-.025em b}\kern-.08em
		T\kern-.1667em\lower.7ex\hbox{E}\kern-.125emX}}
\begin{document}
	\title{Consensus-based Distributed Optimization Enhanced by Integral Feedback}
	\author{Xuan Wang, Shaoshuai Mou, and Brian. D. O. Anderson \vspace{-1.5em}
		\thanks{X. Wang is with the Department of Electrical and Computer Engineering, George Mason University, Fairfax, VA 22030 USA, {xwang64@gmu.edu}. S. Mou is with the School of Aeronautics and Astronautics, Purdue University, West Lafayette, IN 47906 USA {mous@purdue.edu}. His work is supported in part by the NASA University Leadership Initiative (ULI) under grant number 80NSSC20M0161 and Northrop Grumman Corporation.  B. D. O. Anderson is with The Australian National University, Acton, ACT 2601, Australia and Hangzhou Dianzi University, Hangzhou, China, {brian.anderson@anu.edu.au}; his work is supported by Data61-CSIRO and Australian Research Council, Grant  DP190100887. This work has been submitted to the IEEE for possible publication. Copyright may be transferred without notice, after which this version may no longer be accessible.}}
	
\maketitle

\begin{abstract}
	Inspired and underpinned by the idea of integral feedback, a distributed constant gain algorithm is proposed for multi-agent networks to solve convex optimization problems with local linear constraints. Assuming agent interactions are modeled by an undirected graph, the algorithm is capable of achieving the optimum solution with an exponential convergence rate. Furthermore, inherited from the beneficial integral feedback, the proposed algorithm has attractive requirements on communication bandwidth and good robustness against disturbance. Both analytical proof and numerical simulations are provided to validate the effectiveness of the proposed distributed algorithms in solving constrained optimization problems.
\end{abstract}

\begin{IEEEkeywords}
	Distributed Optimization; Integral Feedback; Multi-Agent Networks.
\end{IEEEkeywords}

\section{Introduction}
\IEEEPARstart{C}{ollective} behaviors in nature have motivated rapidly expanding research efforts in the control of multi-agent systems\cite{FJS09DC}. A multi-agent system is composed of multiple interacting subsystems (agents), which makes them more challenging to control than single monolithic systems.
Specifically, the network constraints, stemming from the relations between agents involving sensing, communication or control, usually prohibit the application of traditional methods from controlling the multi-agent systems in a centralized manner.
In order to seek new control approaches that respect the network nature of multi-agent systems, 
\textit{distributed control} has recently received a significant amount of research attention, the goal of which is to allow multi-agent systems to accomplish global objectives through only local coordination. Here, the word `local' connotes interaction between any given agent and a limited number of associated `neighbor' agents, often physically adjacent. 

One of the key problems in multi-agent control is distributed optimization, where  each agent privately processes one objective function and one constraint, and the goal is to minimize the sum of local objective functions \cite{GN17TCNS,TXJ19ARC,KFS18NIPS,XJSM19TAC} subject to all local constraints \cite{PSJW19ARC,KGL19Auto,SJA15TAC,ZSL16Auto,XPY16TAC,QJ15TAC}.
To solve this problem, one applicable approach is the alternating direction method of multipliers (ADMM). While the ADMM method originates from Lagrange duality, it usually needs a centralized state to coordinate across agents \cite{SN11ADMM}. This limits the applicability of the method to fully distributed network systems.
In order to remove the requirement for such a  centralized state, many efforts have been made. For example \cite{GFL16TSP,T16TSP} show the duality can instead be established by introducing slack states on the edges of the network; and in  \cite{T16TSP,NZT17TAC}, the centralized state can be further decomposed by the primal decomposition technique.
Apart from the ADMM methods, another family of distributed coordination method arises from the idea of \textit{consensus}\cite{MAB08JCO}. The aim of consensus is to drive all agents in the network to reach an agreement regarding a certain quantity, which has served as a basis in deriving many distributed algorithms for multi-agent systems such as motion synchronization \cite{FJS09DC}; multi-robot path planning/formation control \cite{XMT15TCNS}; flocking of mobile robots \cite{CH19CSL};  and cooperative sensing \cite{MM10DC}. 
Fitting in the scope of distributed optimization, \textit{consensus} is usually incorporated with \textit{gradient descent} and \textit{projection} operators, which handle objective functions and local constraints, respectively. However, note that in multi-agent systems, different agents may have different local objective functions/constraints.
This means that in general, the consensus, gradient descent, and projection operators will have different equilibria, so that the states in all
agents can never converge to a same point.
To circumvent this difficulty, the work in \cite{AAP10TAC,PWY16Auto,AA09TAC,GBU17TAC} applies a diminishing gain (i.e. $1/t$) to the gradient term in the update equations. As a side effect, this time-variant gain must be shared by all the agents in the network, and the convergence rate of the algorithm will be degraded most commonly to $\mathcal{O}(1/\sqrt{t})$. 
In order to improve the convergence rate, many recent works have shown that a possible approach is to double the dimension of the state vector. See for example the continuous-time update introduced in \cite{BJ14TAC,ZSL16Auto,XPY16TAC,QJ15TAC}, where the role of the extra vector is played by the Lagrangian dual vector for consensus errors;  see also the discrete-time update introduced in  \cite{YAG19Arxiv,WQG15TSP}, where the extra vector performs  gradient tracking;  and the discrete-time update introduced in  \cite{SKA19NIPS,SEK20TAC}, where the extra vectors combines proximal gradient and gradient tracking together.
In these algorithms, the extra states have to be exchanged across the network, necessitating  duplication of the network bandwidth requirement. 
Focusing on the convergence rate and the states to be exchanged across the network, a detailed comparison between this paper and the related works will be provided later in Table \ref{TB_Comp}.





In this paper, we propose a distributed algorithm for constrained optimization that is neither based on diminishing gains nor a doubled dimension of the vectors shared between agents. Actually, by comparing the very fundamental mechanisms of these algorithms, we notice that the common reason why the latter category of algorithms can achieve an improved convergence performance arises from elevating the type of the update to second order, and thereby effectively eliminating the accumulated consensus error.
Inspired by this, in this paper, we propose a continuous-time consensus-based algorithm for constrained distributed optimization based on integral feedback within each agent's controller, and with the integrated signal not being shared with other agents. 
The contributions of the paper in more detail are as follows:
(\textbf{i}) Without a time-variant gain that needs to be shared by agents, the algorithm is capable of achieving { the optimum solution with an asymptotic convergence rate for general convex objective functions with non-unique minimizers}, and a global exponential convergence rate for strongly convex functions.
(\textbf{ii})  {  To achieve exponential convergence, the proposed result only requires \textit{the sum of} all objective functions to be strongly convex at the optimum point\footnote{Please refer to the Assumption \ref{ass2} of the main text for details.}; this is a more relaxed condition compared with the existing results that require all agents' local objective functions to be strongly convex\cite{SJS15Auto,WQG15TSP,SKA19NIPS,ZWM19TSP,SEK20TAC,DDJJ20PIEE}.}
(\textbf{iii}) Inherited from the benefit of integral feedback, it can be theoretically guaranteed  that the proposed algorithm has good robustness against disturbances.
(\textbf{iv}) Apart from requiring the storing at each agent of the integral of the state vector in addition to the state vector itself, { the algorithm does not introduce an extra state vector which has to be exchanged among the agents of the network.} This further distinguishes the work from the existing results based on (primal–dual) saddle point dynamics in \cite{ZSL16Auto,XPY16TAC,QJ15TAC}.
{ Note that the algorithm of this paper is evidently related to the discrete-time algorithm for unconstrained optimization \cite{WQG15TSP,SKA19NIPS,SEK20TAC,DDJJ20PIEE,HZY20Arxiv}.} These authors increased the state dimension by using a form of gradient descent including the last two iterates. On the other hand as noted already, our algorithm is motivated by the very old principle of using integral feedback to cancel steady state errors.  { Different from \cite{SJS15Auto,WQG15TSP,ZWM19TSP,DDJJ20PIEE,HZY20Arxiv}, our algorithm can additionally handle local linear constraints, which commonly exist in many engineering application\cite{J60SIAM,JK11TAC}.}

The rest of the paper is organized as follows. In Section II, we describe the  information flow of the multi-agent network and formulate the problem of constrained distributed optimization. By introducing the  idea of integral feedback, in Section III, we propose a continuous-time  algorithm that can solve distributed optimization problems with linear constraints. The effectiveness, exponential convergence and robustness of the algorithm are theoretically proved in Section IV. Section V provides the numerical validation for the convergence rate and the robustness against disturbance. We finally conclude the paper in Section VI.

\noindent{ \emph{Notation}:} 
Let ${\bf 1}_r$ denote the vector in $\mathbb{R}^r$ with all entries equal to $1$. Let $I_r$ denote the $r\times r$ identity matrix.  We let $\col\{A_1,A_2,\cdots,A_r\}$ be a stack of matrices $A_i$ possessing the same number of columns with the index in a top-down ascending order, $i=1,2,\cdots,r$.
By $M^{\top}$ is meant the transpose of a matrix $M$. Let $\ker M$ and $\image M$ denote the kernel  and image of a matrix $M$, respectively. Let $\otimes$ denote the Kronecker product.

\begin{table*}[b]
	\caption{Comparison with existing algorithms}\label{TB_Comp}
	\vspace{-2em}
	\begin{center}
		\begin{tabular}[t]{c|c|c|c|c|c}
			\hline
			\hline
			\textbf{Algorithm} & \textbf{Key idea}  & \textbf{State to share}  &\textbf{Local constraints} &\textbf{Update form} & \textbf{Exponential Convergence}\\
			\hline
			The proposed algorithm & Integral Feedback & Dimension $n$& Linear, closed & Continuous &  Yes\\
			\hline
			{Algorithms in \cite{PWY16Auto} } & Diminishing step-size/gain  & Dimension $n$ & Compact& Discrete & No, $\mathcal{O}(1/\sqrt{t})$ \\ \hline
			{ Algorithms in \cite{SJS15Auto}} & { Integral Feedback$^{\dagger}$}  & { Dimension $n$} & { Not applicable}& { Continuous} & { Yes}\\\hline
			{ Algorithms in \cite{ZWM19TSP,WQG15TSP,DDJJ20PIEE} }& { Gradient tracking} & { Dimension $2n^\star$} & { Not applicable}& { Discrete} & { Yes}\\\hline
			Algorithm in \cite{BJ14TAC} & Saddle point dynamics & Dimension $2n$ & Not applicable& Continuous & No theoretical guarantee\\\hline
			Algorithms in \cite{ZSL16Auto,XPY16TAC,QJ15TAC} & Saddle point dynamics & Dimension $2n$ & Closed& Continuous & No theoretical guarantee\\\hline
		\end{tabular}
	\end{center}
	\hspace{2em} \scriptsize 
	{ $\dagger$: Note that \cite{SJS15Auto} requires an extra process for designing control gains. In this paper, the process is not required.}
	
	{ \hspace{2em} \scriptsize  ~$\star$: In \cite{WQG15TSP}, certain choices for its parameter design can lead to a discrete counterpart of the result in this paper. In this case, its \textit{state to share} is reduced to $n$.}
\end{table*}
\vspace{-0.5em}

\section{Problem Formulation}

Consider a network of ${m}$ agents in which each agent $i$ is able to communicate with certain other nearby agents called its neighbors, denoted by $\mathcal{N}_i$. The neighbor relations can be described by a graph $\mathbb{G}$ such that there is an edge from $j$ to $i$ if and only if $j\in \mathcal{N}_i$.
We assume $\mathbb{G}$ is \textit{connected} and \textit{undirected}.
Associated with each agent is a local state $x_i\in\mathbb{R}^n$; a convex function $f_i(\cdot):\mathbb R^n\rightarrow\mathbb R$; and a linear constraint $A_ix_i=b_i$, where $A_i\in\mathbb{R}^{n_i\times n}$ and $b_i\in \image A_i \subset\mathbb{R}^{n_i}$. The problem of \textit{interest} is to develop a \textit{distributed} algorithm which enables all nodes of $\mathbb{G}$ to reach a consensus value solving the problem\vspace{-0.2em}
\begin{align}\label{eqobjective}
\text{minimize} \qquad&\sum_{i=1}^{m} f_i(x_i).\\
\text{subject to} \qquad &A_i x_i=b_i, ~~i=1\cdots,m \label{eq_constraint}\\
 &x_1=x_2=\cdots=x_m. \label{eq_consensus}
\end{align}
\begin{remark}\label{Rm_LinearCons}
The linear constraint in \eqref{eq_constraint} arises naturally from many engineering applications\cite{J60SIAM,JK11TAC}. Here, to avoid trivialities, we assume $b_i\in \image A_i$, and  $\rank(A)<n$, where $A=\col\{A_1,\cdots,A_m\}$. This guarantees the optimization domain defined by equations \eqref{eq_constraint}-\eqref{eq_consensus} is non-empty and non-unique.   
\end{remark}


\vspace{-0.5em}
\section{The Algorithm and Main Result}
 In order to solve the consensus-based constrained optimization problem characterized by \eqref{eqobjective}-\eqref{eq_consensus}, we first introduce some preliminary results.\vspace{-1em}
\subsection{Preliminary Results}
\noindent\textbf{Distributed Consensus:} Consider a continuous-time update of the following form \cite{L04CDC}
\begin{align}\label{eq_Alg_consensus}
\dot{x}_i=-\sum_{j\in\mathcal{N}_i} (x_i -x_j),
\end{align}
where each agent in the network tries to reduce the distances between itself and its neighbors.  Update \eqref{eq_Alg_consensus} will drive all states to a consensus value, i.e. there exists a certain $\bm{x}^*=\bm{1}_m\otimes x^*$, such that $x_i-x^*\rightarrow 0$ for all $i$ and the convergence is exponentially fast. 

\noindent\textbf{Consensus-based distributed optimization:}
Leaving aside temporarily the constraints \eqref{eq_constraint}, the  attempt to find an algorithm also achieving the optimization objective \eqref{eqobjective} requires the introduction of a gradient term, so that the update of each agent becomes\vspace{-0.2em}
\begin{align}\label{eq_Alg_opt}
\dot{x}_i= -\alpha(t)\nabla f_i(x_i)-\sum_{j\in\mathcal{N}_i} {s_{ij}} (x_i -x_j),
\end{align}
where $\alpha(t)$ is a positive gain shared by all the agents and $\nabla f_i(x_i)$ is the gradient\footnote{Note that if the function is not differentiable, the gradient can usually be replaced by a sub-gradient of the function with the convergence property being unchanged. However, since sub-gradients are not unique, this can lead to non-unique trajectories and requires the tool of Filippov-solution for analysis.
 of $f_i(x_i)$, $s_{ij}>0$ is a parameter, which is usually associated with the edge weight of the network. 
Note that in \eqref{eq_Alg_opt}, if $\alpha(t)$ is chosen as a fixed positive constant, unless all $f_i$ are minimized by a common vector, there may not exist a  steady state $x^*$ for equation \eqref{eq_Alg_opt} which also follows the consensus property \eqref{eqobjective}.}. 

In order to solve the consensus-based optimization problem, a novel algorithm is presented in \cite{AA09TAC,AAP10TAC,PWY16Auto}, where the authors apply a diminishing gain $\alpha(t)$ to the discrete-time version of \eqref{eq_Alg_opt} to eliminate the consensus error. Furthermore, with introduction of  an additional projection operator, this algorithm is also able to handle local constraints. 
A continuous version of the algorithm developed in \cite{AA09TAC,AAP10TAC,PWY16Auto} is \vspace{-0.2em}
\begin{align}\label{eq_Alg_opt_cons_dis}
\dot{x}_i= \mathcal{P}_i\left[-\alpha(t)\nabla f_i(x_i)-\sum_{j\in\mathcal{N}_i} (x_i -x_j)\right],
\end{align}
where $\mathcal{P}_i[\cdot]$ is a projection operator that projects any vector to the tangent space\footnote{For the linear constraints, one has $\mathcal{P}_i[s]=P_i\cdot s$, where $P_i\in\mathbb{R}^{n\times n}$ is a projection matrix to $\ker A_i$.} of the agent's local constraint at the point $x_i(t)$, and it guarantees  that ${x}_i$ always satisfies the local constraint of agent $i$. By letting all agents  share a diminishing gain such that $\alpha(t)\to 0$ and $\int_0^{\infty}\alpha(t)\to\infty$, it has been theoretically proved that the states  will asymptotically reach a consensus at the minimizer of $\sum_{i=1}^m f_i(x)$ subject to all agents' local constraints\cite{AAP10TAC,PWY16Auto}. Meanwhile, since the effect of the gradient term is discounted by the diminishing gain, the convergence rate of the algorithm is at most $\mathcal{O}(1/\sqrt{t})$. Here, to improve the convergence rate, the key idea is to get rid of the diminishing gain $\alpha(t)$ in \eqref{eq_Alg_opt_cons_dis}. This has led to the main result of this paper.

\vspace{-0.5em}
\subsection{The proposed update and  main result}
In this paper, instead of using a diminishing gain or an additional state vector to each agent which has to be exchanged with neighbors, our key idea stems from  introducing an additional integral term to effectively eliminate the steady state error on consensus. We propose the following continuous-time distributed algorithm,\vspace{-0.3em}
\begin{align}\label{Dupdatex}
\dot x_i= -P_i\left(\nabla f_i(x_i)+\sum_{j\in\mathcal{N}_i}(x_i-x_j)+ \int_{0}^{t}\sum_{j\in\mathcal{N}_i}(x_i-x_j)\right) 
\end{align}
where $x_i(0)$ are initialized such that $A_ix_i(0)=b_i$; and $P_i\in\mathbb{R}^{n\times n}$ is a projection matrix to $\ker A_i$.  
\begin{remark}
	Obviously, the proposed algorithm  is distributed, because the state update of each agent only relies on the information of itself and that of its neighbors.
	In update \eqref{Dupdatex}, the projection matrix $P_i$ and the special initialization on $x_i(0)$ are used to handle the linear constraint \eqref{eq_constraint}. As a special case, if for one or more agents, the linear constraint does not exist, then one can correspondingly initialize $x_i(0)$ as an arbitrary value and replace the  projection matrix by an identity matrix.  Further note that similarly to \cite{SJS15Auto,A17SIAM}, the integral term requires extra historical state information; as will be seen later in equation \eqref{Cupdatey}, such information can be equivalently stored in a local state $y_i=\int_{0}^{t}\sum_{j\in\mathcal{N}_i}(x_i-x_j)$, and moreover, this extra state does not need to be shared across the network.
\end{remark}

Ahead of studying the properties of Algorithm \eqref{Dupdatex}, we introduce the following assumptions.
\begin{assumption}\label{ass1} 
	For all $i=1,\cdots,m$, suppose $f_i(\cdot)$ is convex, continuously differentiable and its gradient is Lipschitz.
	Suppose $\rank(A)<n$, where $A=\col\{A_1,\cdots,A_m\}$, such that the feasible set defined by the constraints \eqref{eq_constraint}-\eqref{eq_consensus} is non-empty. Suppose a solution $x^*$ to \eqref{eqobjective}-\eqref{eq_consensus} exists, where $x^*$ may be non-unique\footnote{  We acknowledge the suggestion of an anonymous reviewer to include a non-uniqueness condition in this assumption.}.
\end{assumption}
\begin{assumption}\label{ass2}
	Suppose Assumption \ref{ass1} holds. Further suppose $F(x)=\sum_{i=1}^m f_i(x)$ is strongly convex\footnote{A function ${\displaystyle F(\cdot)}$  is called strongly convex at $x^*$ with parameter $\omega > 0$ if $a^{\top}[\nabla F(a+x^*)-\nabla F(x^*)]\geq \omega\|a\|_2^2$ holds for all vectors $a$ in its domain.} at $x^*$, where $x^*$ is the unique minimizer to problem \eqref{eqobjective}-\eqref{eq_consensus}.
\end{assumption}

\smallskip
The Algorithm \eqref{Dupdatex}, along with assumptions \ref{ass1} and \ref{ass2} allows us to propose the following Theorem. 

\smallskip

\begin{theorem}\label{thm1}
	Suppose the network $\mathbb{G}$ is connected and undirected; suppose Assumption \ref{ass1} holds.
	Then given any $x_i(0)$ such that $A_ix_i(0)=b_i$, update \eqref{Dupdatex} drives the states $x_i(t)$ of all agents asymptotically { to a point $x^*$, which is an optimum solution to \eqref{eqobjective}-\eqref{eq_consensus}.}  Furthermore, if Assumption \ref{ass2} also holds, the convergence is globally exponentially fast.
\end{theorem}

\smallskip
\begin{remark}\label{RThm}
	Note immediately that in Assumption \ref{ass2}, we only require the sum $F(x)=\sum_{i=1}^m f_i(x)$ to be strongly convex (at $x^*$), which is a more relaxed condition compared with most of the existing results\cite{WQG15TSP,SKA19NIPS,SEK20TAC} that require each single $f_i(x_i)$ to be strongly convex. This allows the algorithm to handle the scenarios when linear objective functions or exponential penalty functions are associated with some of the agents.
\end{remark}

\begin{corollary}\label{Cl1}
	Suppose Assumption \ref{ass2}  holds. Consider the update under disturbance,\vspace{-0.3em}
	\begin{align}\label{DupdatexD}
		\dot x_i\!=\! -P_i\!\left(\!\nabla f_i(x_i)\!+\!\!\sum_{j\in\mathcal{N}_i}\!(x_i\!-\!x_j)\!+\!\! \int_{0}^{t}\!\sum_{j\in\mathcal{N}_i}\!(x_i\!-\!x_j)\!\right) +\!v_i(t)\!
	\end{align}
where $v_i(t)$ is a disturbance that is bounded for all $i=1,\cdots,m$ and $t\ge0$. Then, the proposed algorithm \eqref{Dupdatex} is robust against bounded disturbance, that is, as $t\to\infty$, one has that $\|x_i(t)-x^*\|$ remains bounded for all $i=1,\cdots,m$. 
\end{corollary}

\begin{remark}\label{RCol}
	In Corollary \ref{Cl1}, the disturbance $v_i(t)$ can frequently be caused by communication issues or sensor mismatches\cite{TM12ACC}. 	
	Also note that the robustness property proposed here is primarily a consequence of the exponential convergence of algorithm \eqref{Dupdatex}. On the contrary, for distributed algorithms that only have asymptotic convergence (i.e. $\mathcal{O}(1/t)$, $\mathcal{O}(1/\sqrt{t})$), the same disturbance will lead the output error to be unbounded. This statement will be validated later by simulations.
\end{remark}

\vspace{-0.5em}
\section{Analysis}
This section proves the main results of the paper. \vspace{-0.5em}
\subsection{Steady-state Analysis}
{ We first propose the following lemma, which shows that the equilibrium point of \eqref{Dupdatex} exists, and it is consistent with the solution to problem \eqref{eqobjective}-\eqref{eq_consensus}.}
\begin{lemma}\label{lm1}	
	Consider the updates \eqref{Dupdatex} under Assumption \ref{ass1}, where $\mathbb{G}$ is connected and undirected.
	Then an equilibrium point $x_i^*$ of equation \eqref{Dupdatex} exists. Furthermore, for all $i=1,\cdots,m$, there holds  $x_i^*=x^*$, ensuring that the equilibrium point obeys the consensus property and optimizes the constrained optimization problem \eqref{eqobjective}-\eqref{eq_consensus}. 
\end{lemma}
\smallskip

For simplicity in analyzing the proposed algorithm from a global perspective, define $\bm{x}=\col\{x_1, \cdots, x_m \}\in\mathbb{R}^{mn}$, $\bar{P}=\diag\{P_1,\cdots,P_m\}\in\mathbb{R}^{mn\times mn}$, {$\nabla f(\bm{x})=\col\{\nabla f_1(x_1),\cdots,\nabla f_m(x_m)\}$}, and $\bar{L}=L \otimes I_n \in\mathbb{R}^{mn\times mn}$, where $L\in\mathbb{R}^{m\times m}$ is the Laplacian matrix of the graph $\mathbb{G}$\footnote{Since $\mathbb{G}$ is connected and undirected, $L$ must be symmetric with kernel spanned by ${\bf{1}_m}$, the $m$-vector of all 1's.}. Then update \eqref{Dupdatex} can be rewritten as 
\begin{align}
\dot{\bm{x}}&= -\bar{P}\left(\nabla f(\bm{x})+\bar{L}\bm{x}+ \int_{0}^{t}\bar{L}\bm{x}\right).\label{Iupdatex}
\end{align}
This is further equivalent to
\begin{align}
\dot{\bm{x}}&= -\bar{P}\left(\nabla f(\bm{x})+ \bar{L}\bm{x}+ \bm{y}\right)\label{Cupdatex}\\
\dot{\bm{y}}&= \bar{L}\bm{x}\label{Cupdatey}
\end{align}
where $\bm{y}\in\mathbb{R}^{mn}$ and  $\bm{y}(0)=0$.
\begin{remark}
	Evidently the dynamics \eqref{Cupdatex}-\eqref{Cupdatey}  is equivalent to \eqref{Dupdatex}, where the integral term implicitly introduces an extra state $\bm{y}$. Each component of this extra state can be obtained via local computations and stored  by each agent,  and does not have to be exchanged across the network. 
	For existing algorithms characterized by saddle-point dynamics, the extra states must be exchanged across the network \cite{BJ14TAC,AAW17JOpt}.   
\end{remark}

\smallskip

\noindent{\bf{Proof of Lemma \ref{lm1}:}} 
To prove Lemma \ref{lm1}, since (\ref{Cupdatex})-(\ref{Cupdatey}) and \eqref{Dupdatex} are equivalent, it is sufficient to show that there exist equilibrium points $(\bm{x}^*,\bm{y}^*)$ such that 
\begin{align}
0&= -\bar{P}\left(\nabla f(\bm{x}^*)+\bar{L}\bm{x}^*+ \bm{y}^*\right)\label{eq_eqlix}\\
0&= \bar{L}\bm{x}^* \label{eq_eqliy}
\end{align}
where $\bm{x}^*=\bm{1}_m\otimes x^*$ and $x^*$ is a minimizer of $F(x)=\sum_{i=1}^m f_i(x)$ subject to $A_ix^*=b_i$ for all $i$; $\nabla f(\bm{x}^*)$ is the column stack of the vectors $\left.\frac{d f_i(x_i)}{d x_i}\right|_{x^*}$.

\smallskip
\noindent \textbf{Existence}: 
	Based on Assumption \eqref{ass1}, let $x^*$ be a minimizer of $\sum_if_i(x)$ subject to $A_ix=b_i$, $i=1\cdots,m$. Then, by standard Lagrange multiplier theory\cite{D97Nonlinear}, there exist Lagrange multipliers $z_i^*$  such that 
	\begin{equation}\label{eq_dfiaz}
	\sum_{i=1}^m\nabla f_i(x^*)-\sum_{i=1}^mA_i^{\top}z_i^*=0
	\end{equation}
	Using the $z_i^*$, we make the definition
	\begin{equation}\label{eq_yi*}
	y_i^*=A_i^{\top} z_i^*-\nabla f_i(x^*).   
	\end{equation}	
Based on $x^*$ and $y_i^*$, let
\begin{align}\label{eq_defelu}
	\bm{x}^*=\bm{1}_m\otimes x^*\quad \text{ and } \quad\bm{y}^*= \col\{y_1^*,\cdots,y_m^*\}.
\end{align}
In the following, we prove that the equilibrium to dynamics \eqref{Cupdatex}-\eqref{Cupdatey} exists, by showing that $\bm{x}^*$ and $\bm{y}^*$ satisfy equations \eqref{eq_eqlix}-\eqref{eq_eqliy}. First, since $\ker \bar{L}= \image (\bm{1}\otimes I_n)$ and $\bm{x}^*=\bm{1}_m\otimes x^*$, one has $\bar{L}\bm{x}^*=0$, which is equation \eqref{eq_eqliy}. 

To continue, from equation \eqref{eq_yi*}, and the definitions of $\bm{x}^*$ and $\bm{y}^*$ in \eqref{eq_defelu}, one has 
	\begin{equation}\label{eq_fbA}
	\nabla f(\bm{x}^*)+ \bm{y}^*=\bar A^{\top}{\bm{z}^*}
	\end{equation}
	where $\bar{A}=\diag\{A_1,\cdots, A_m\}$. Indeed
	\begin{equation}
	\nabla f(\bm{x}^*)+\bar L{\bm{x}^*}+ \bm{y}^*=\bar A^{\top}{\bm{z}^*} 
	\end{equation}	
	From this, since $\bar P\bar A^{\top}=0$,  we obtain
	\begin{equation}\label{eq_exiseq}
	0=-\bar P(\nabla f(\bm{x}^*)+\bar L\bm{x^*}+ \bm{y}^*)
	\end{equation}
	This ensures the  satisfaction of \eqref{eq_eqlix}. Furthermore, since $\bm{y}(0)=0$, the integration of $\bm{y}$ in equation \eqref{Cupdatey} leads to another implicit condition, that is for all $t$, $\bm{y}(t)\in\image \bar{L}$. Hence, $\bm{y}^*\in\image \bar{L}$. To  validate this, from \eqref{eq_dfiaz} and \eqref{eq_yi*} one has 
	\begin{align*}
	\sum_{i=1}^m y_i^*&=\sum_{i=1}^m \left[ A_i^{\top} z_i^*-\nabla f_i(x^*)\right]\\
	&=-\left(\sum_{i=1}^m\nabla f_i(x^*)-\sum_{i=1}^mA_i^{\top}z_i^*\right)=0
	\end{align*}
	This ensures the  satisfaction of $\bm{y}^* \in \mbox{image }\bar L$ and establishes the existence of the equilibrium points $(\bm{x}^*,\bm{y}^*)$.

	\smallskip
	\noindent { \textbf{Consistency}}:
	Based on the existence of the equilibrium $(\bm{x}^*,\bm{y}^*)$, we now show that any $\bm{x}^*$ satisfying \eqref{eq_eqlix}-\eqref{eq_eqliy} takes the form set out in the lemma statement and is consistent with the solution to problem \eqref{eqobjective}-\eqref{eq_consensus}.	
	{  Since $\ker \bar{L}= \image (\bm{1}\otimes I_n)$ and $\bar{L}\bm{x}^*=0$, one has $\bm{x}^*=\bm{1}_m\otimes u$, where $ u\in\mathbb{R}^{n}$. This directly leads to the consensus property. To continue, we only need to show that $u={x}^*$ is a solution to problem \eqref{eqobjective}-\eqref{eq_consensus}.
	}	
	To do this, recall that $\image P_i=\ker A_i$, thus, $\ker P_i=\image A_i^{\top}$. From equation \eqref{eq_eqlix} and 
	$\bar{P}=\diag\{P_1,\cdots,P_m\}$, there exists a $\bm{z}^*= \col\{z_1^*,\cdots,z_m^*\}$, $z^*_i\in\mathbb{R}^{n_i}$ such that 
	\begin{align}\label{eq_Az}
	\nabla f(\bm{x}^*)+\bar{L}\bm{x}^*+ \bm{y}^*=\bar{A}^{\top}\bm{z}^*.
	\end{align}
	Further recall that $\nabla f(\bm{x}^*)=\nabla f(\bm{1}_m\otimes u)$ is the column stack of the vectors $\frac{d f_i(x_i)}{d x_i}|_{u}$, then multiplying equation \eqref{eq_Az} on the left by $(\bm{1}_m\otimes I_n)^{\top}$  yields
	\begin{equation}\label{eq_sumDfx}
	(\bm{1}_m\otimes I_n)^{\top}\nabla f(\bm{x}^*)=\sum_{i=1}^m\nabla f_i(u)=\sum_{i=1}^m A_i^{\top}z_i.
	\end{equation}
	This, by standard Lagrange multiplier theory, tells us  $u$  is a critical point for $F(x)=\sum_{i=1}^mf_i(x)$ on the manifold defined by $A_ix=b_i$ for all $i$. 
	Thus, $u=x^*$ is a minimizer to problem \eqref{eqobjective}-\eqref{eq_consensus}.
	This completes the proof. \qed

	{ As a side remark, we observe that under Assumption \ref{ass1}, a solution $x^*$ can be non-unique, but under Assumption \ref{ass2}, ${x}^*$ is clearly unique. }
	Further more, note that the uniqueness, or otherwise, of $y_i^*$ does not influence the result of Lemma \ref{lm1}, as it does not originally appear in the update \eqref{Dupdatex}. Actually in \eqref{eq_dfiaz}, if $\begin{bmatrix}
	A_1^{\top}, \cdots, A_m^{\top}
	\end{bmatrix}$
	does not have linearly independent columns, the value of $z_i^*$ is non-unique \cite{D97Nonlinear}. Consequently, the value of $y_i^*$ is non-unique.

\subsection{Change of coordinate frame}
In order to examine the transient behavior of (\ref{Cupdatex})-(\ref{Cupdatey}), it is convenient to make a coordinate transformation which ensures that in the new coordinates, the equilibrium point corresponding to $\bm{x}^*$ is moved to zero. Here, we change the origin of updates (\ref{Cupdatex})-(\ref{Cupdatey}), by defining vectors $\tilde{\bm{x}}, \tilde{\bm{y}}$ as \vspace{-0.3em}
\begin{align}\label{eq_deftdxy}
\tilde{\bm{x}}=\bm{x}-\bm{x}^*\nonumber\\
\tilde{\bm{y}}=\bm{y}-\bm{y}^*
\end{align}
where, as above, $\bm{x}^*=\bm{1}_m\otimes x^*$ and $\bm{y}^*$ satisfies \eqref{eq_exiseq} with a certain Lagrange multiplier $\bm{z}^*$. Note that when $\bm{x}^*$ and $\bm{y}^*$ are non-unique, one can make an arbitrary choice consistent with \eqref{eq_eqlix}-\eqref{eq_eqliy}. Evidently, 
\begin{align}
\dot{\tilde{\bm{x}}}=&-\bar P\left(\nabla f(\tilde{\bm{x}}+\bm{x}^*)-\bar{L}(\tilde{\bm{x}} +\bm{x}^*)- (\tilde{\bm{y}}+\bm{y}^*)\right)\nonumber\\
=&-\bar P\left(\left[\nabla f(\tilde{\bm{x}}+\bm{x}^*)-\nabla f(\bm{x}^*)\right]-\bar{L}\tilde{\bm{x}}- \tilde{\bm{y}}\right)\label{eq:newcordx}\\
\dot{\tilde{\bm{y}}}=&~~\bar{L}(\tilde{\bm{x}} +\bm{x}^*)
=\bar{L}\tilde{\bm{x}}\label{eq:newcordy}
\end{align}

\smallskip
To continue, we further modify updates (\ref{eq:newcordx})-(\ref{eq:newcordy}) by a frame transformation. Since the linear equation set $A_ix=b_i$, $i=1,\cdots,m$ has multiple solutions,  there exists {at least one nonzero vector} in the kernel of every $A_i$, i.e. in the range of every $P_i$. Observe that if $u$  is such a vector, then $P_iu=u$, $\forall i \in\{1,\cdots,m\}$. Recalling that $\ker \bar{L}=\image (\bm{1}\otimes I_n)$, then, $\bar{P}\bar{L}\bar{P}(\bm{1}\otimes u)=0$, which means $\bar{P}\bar{L}\bar{P}$ is singular. Thus, there exists an orthogonal matrix $Q=\begin{bmatrix}
R_1&R_2
\end{bmatrix}$, with $R_1\in\mathbb{R}^{mn\times\bar{n}_1}$, $R_2\in\mathbb{R}^{mn\times\bar{n}_2}$, $\bar{n}_1+\bar{n}_2=mn$, such  that 
\begin{align}\label{eq_rotation}
Q^{\top}\bar{P}\bar{L}\bar{P}Q=\begin{bmatrix}
0&&0\\
0&&R_2^{\top}\bar{P}\bar{L}\bar{P}R_2
\end{bmatrix},
\end{align}
where the matrix $R_2^{\top}\bar{P}\bar{L}\bar{P}R_2$ is non-singular. Now define new vectors $X,Y$ by the transformations
\begin{equation}\label{eq_defXY}
X=Q^{\top}\tilde{\bm{x}} ,\;\;Y=Q^{\top}\bar{P}\tilde{\bm{y}}.
\end{equation}
Multiplying the differential equations  (\ref{eq:newcordx})-(\ref{eq:newcordy}) on the left, by $Q^{\top}$ and $Q^{\top}\bar P$, respectively, yields
\begin{align}\label{eq:XYde}
\dot X=&- Q^{\top}\bar P[\nabla f(QX+\bm{x}^*)-\nabla f(\bm{x}^*)]-Q^{\top}\bar P\bar{L}\bar PQX\nonumber\\
&- Y
\\\label{eq:XYde2}
\dot Y=&~~Q^{\top}\bar P\bar{L}\bar PQX
\end{align}
Note that in the derivation of \eqref{eq:XYde}-\eqref{eq:XYde2}, we have replaced $QX$ with $\bar{P}QX$. This equality holds because both $x_i(t)$ for all $t$ and $x^*$ are solutions to $A_ix=b_i$; then $P_i(x_1-x^*)=(x_1-x^*)$, that is, $\bar{P}QX=\bar{P}\tilde{\bm{x}}=\col\{P_1(x_1-x^*),\cdots,P_m(x_m-x^*)\}=\tilde{\bm{x}}=QX$. Based on \eqref{eq:XYde}-\eqref{eq:XYde2}, further partition the vectors $X,Y$ as
\begin{equation}
X=\begin{bmatrix}
X_1\\X_2
\end{bmatrix},\quad Y=\left[\begin{array}{c}Y_1\\Y_2\end{array}\right]
\end{equation}
where $X_1,Y_1\in\mathbb{R}^{\bar{n}_1}$ and $X_2, Y_2 \in \mathbb R^{\bar{n}_2}$.
Consider now the equations for $X_1,Y_1$ alone. Using equation \eqref{eq_rotation}, there results
\begin{align}\label{eq_X_1}
\dot X_1=&- R_1^{\top}\bar P\left([\nabla f(\tilde{\bm{x}}+\bm{x}^*)-\nabla f(\bm{x}^*)]- Y_1\right)\\
\dot Y_1=&~~0 \label{eq_Y_1}
\end{align}
and
\begin{align}\label{eq_X2}
\dot X_2 =&- R_2^{\top}\bar P[\nabla f(QX\!+\!\bm{x}^*)-\nabla f(\bm{x}^*)]\nonumber\\
&-R_2^{\top}\bar{P}\bar{L}\bar{P}R_2X_2- Y_2\\
\dot Y_2 \!=&~~R_2^{\top}\bar{P}\bar{L}\bar{P}R_2 X_2  \label{eq_Y2}
\end{align}
Now observe that $Y_1(0)=0$. The argument is as follows. Because  $Q=\begin{bmatrix}R_1&R_2 \end{bmatrix}$, from equation \eqref{eq_rotation}, one has $R_1^{\top}\bar{P}\bar{L}=0$. Recall that $\bm{y}(0)=0$ and $\bm{y}^*\in\image \bar{L}$, then  $Y_1(0)=R_1^{\top}\bar{P}\tilde{\bm{y}}(0)=R_1^{\top}\bar{P}\left(\bm{y}(0)-\bm{y}^*\right)=0$. In light of \eqref{eq_Y_1}, this means that  $Y_1=0$ for all $t$, furthermore, $X^{\top}Y=X_1^{\top}Y_1+X_2^{\top}Y_2=X_2^{\top}Y_2$.

Stability questions concerning \eqref{eq:newcordx} and \eqref{eq:newcordy} thus can be treated as stability questions concerning the equations \eqref{eq_X_1}, \eqref{eq_X2} and \eqref{eq_Y2}, leaving out \eqref{eq_Y_1}.

\vspace{-0.5em}
\subsection{Proof of Theorem \ref{thm1}}
The proof comprises two main steps. In the first step, under Assumption \ref{ass1}, Lyapunov theory is used to establish asymptotic stability of \eqref{eq_X_1}, \eqref{eq_X2} and \eqref{eq_Y2}. In the second step, the equilibrium is shown to be {globally}  exponentially stable provided Assumption \ref{ass2} holds. 

\noindent\textbf{Asymptotic Stability:} Noting that $R_2^{\top}\bar{P}\bar{L}\bar{P}R_2$ is symmetric  positive definite, we can define a positive definite function $V$ of $X$ and $Y_2$ as follows:
\begin{equation}\label{eq_defV}
V(X,Y_2)=\frac{1}{2}\left[X^{\top}X+ Y_2^{\top}(R_2^{\top}\bar{P}\bar{L}\bar{P}R_2)^{-1}Y_2\right].
\end{equation}
Computing the derivative along motions of equations (\ref{eq:XYde})-(\ref{eq:XYde2}) gives us 
\begin{align}\label{eq_Vdot}
\dot V\!=&- X^{\top}\!Q^{\top}\!\bar P[\nabla \!f(QX\!+\!\bm{x}^*)\!-\!\nabla \!f(\bm{x}^*)]\!-\!X^{\top}\! Q^{\top}\!\bar{P}\bar{L}\bar{P}QX\nonumber\\
&-  X^{\top} Y+ Y_2^{\top}(R_2^{\top}\bar{P}\bar{L}\bar{P}R_2)^{-1}(R_2^{\top}\bar{P}\bar{L}\bar{P}R_2)X_2\nonumber\\
=&- X^{\top}\!Q^{\top}\!\bar P[\nabla \!f(QX\!+\!\bm{x}^*)\!-\!\nabla \!f(\bm{x}^*)]\!-\!X^{\top}\! Q^{\top}\!\bar{P}\bar{L}\bar{P}QX\nonumber\\
=&- X^{\top}\!Q^{\top}\![\nabla\! f(QX\!+\!\bm{x}^*)\!-\!\nabla\! f(\bm{x}^*)]
\!-\!X^{\top}\! Q^{\top}\!\bar{L}QX
\end{align}
Note  that the last equality holds because $QX=\bar{P}QX$ and  $X^{\top}Y=X_1^{\top}Y_1+X_2^{\top}Y_2=X_2^{\top}Y_2$.

Since $\nabla f(\bm{x})=\col\{\nabla f_1(x_1),\cdots,\nabla f_m(x_m)\}$ and each $f_i(\cdot)$ is convex, one has $f(\cdot)$ is also convex, that is $- X^{\top}\!Q^{\top}[\nabla f(QX+\bm{x}^*)-\nabla f(\bm{x}^*)]\le0$. Thus,  $\dot V$ in \eqref{eq_Vdot} is non-positive. By applying LaSalle's  Theorem \cite{H02Nonlinear}, we know the system converges to $\dot V=0$, i.e. $\bar{L}Q{X}=X^{\top}\!Q^{\top}\![\nabla\! f(QX\!+\!\bm{x}^*)-\nabla\! f(\bm{x}^*)]=0$. Because $\ker \bar{L}=\image \left(\bm{1}_m\otimes I_n\right)$, one has $Q{X}=\bm{1}_m\otimes q$, $ q\in\mathbb{R}^{n}$. 
Consequently, 
\begin{align}\label{eq_XQdQX}
	&X^{\top}\!Q^{\top}[\nabla f(QX+\bm{x}^*)-\nabla f(\bm{x}^*)]\nonumber\\
	=~&[\bm{1}_m\otimes q]^{\top}[\nabla f(\bm{1}_m\otimes (q+{x}^*))-\nabla f(\bm{1}_m\otimes {x}^*)]\nonumber\\
	=~&q^{\top}\sum_{i=1}^m\left[\nabla f_i(q+x^*)-\nabla f_i(x^*)\right]\nonumber\\
	=~&{ q^{\top}\left[\nabla F(q+x^*)-\nabla F(x^*)\right]=0}
\end{align}
{ The last equality holds because $F(x)=\sum_{i=1}^m f_i(x)$. Recall that $\bm{x}-\bm{x}^*=Q{X}=\bm{1}_m\otimes q$ and $\bm{y}-\bm{y}^*=\bar{L}Q{X}=0$, thus, the states of all agents converge asymptotically to a same vector $x^*+q$, which satisfies the consensus property \eqref{eq_consensus}. To complete the proof, given $x^*$ is a minimizer to \eqref{eqobjective}-\eqref{eq_consensus}, we only need to show $x^*+q$ is also a minimizer. To do this, because $\bar{P}QX=QX$, then $q\in\cap_{i=1}^{m}\ker(A_i)$. Clearly $x_i=x^*+q$ satisfies constraint \eqref{eq_constraint}. Furthermore, since $x^*$ is a minimizer, there holds $q^{\top}\sum_{i=1}^m\nabla f_i(x^*)=q^{\top}\nabla F(x^*)=0$. This, together with $\eqref{eq_XQdQX}$ yields $q^{\top}\nabla F(q+x^*)=0$. Finally, due to the convexity of $F(x)$, one has $F(x^*)\ge F(x^*+q)+(-q)^{\top}\nabla F(q+x^*)=F(x^*+q)$. Thus, $F(x^*+q)=F(x^*)$ and $x^*+q$ must be a minimizer to problem \eqref{eqobjective}-\eqref{eq_consensus}. This establishes the asymptotic convergence of the algorithm.}	

\smallskip

\noindent\textbf{Exponential Stability:} Here, we  establish  the  exponential  convergence  rate  of the  proposed  update under Assumption \ref{ass2}. { Note that in this case $x^*$ must be unique.} Define 
\begin{align}\label{eq_defTheta}
	\Theta_1= \bar{P}Q\begin{bmatrix}
		X_1\\0_{\bar{n}_2}
	\end{bmatrix}\quad\text{and}\quad \Theta_2= \bar{P}Q\begin{bmatrix}
		0_{\bar{n}_1}\\X_2
	\end{bmatrix}
\end{align}	
such that $\Theta_1, \Theta_2 \in\mathbb{R}^{mn}$ and $\Theta_1+\Theta_2=\bar{P}Q\begin{bmatrix}
	X_1\\X_2
\end{bmatrix}=\bar{P}QX$. To continue, recall the second last line of equation \eqref{eq_Vdot} and the fact that $QX=\bar{P}QX$, then $\Theta_1+\Theta_2=QX$ and
\begin{align}\label{eq_Vdotexp}
\dot V\!=&- X^{\top}\!Q^{\top}\!\bar P[\nabla \!f(\bar PQX\!+\!\bm{x}^*)\!-\!\nabla \!f(\bm{x}^*)]\!-\!X^{\top}\! Q^{\top}\!\bar{P}\bar{L}\bar{P}QX\nonumber\\
=&- (\Theta_1+\Theta_2)^{\top}[\nabla \!f(\Theta_1+\Theta_2+\!\bm{x}^*)\!-\!\nabla \!f(\bm{x}^*)]\nonumber\\
&-\!X^{\top}\! Q^{\top}\!\bar{P}\bar{L}\bar{P}QX
\end{align}
Since $f(\cdot)$ is convex, $-(\Theta_1+\Theta_2)^{\top}[\nabla \!f(\Theta_1+\Theta_2+\!\bm{x}^*)\!-\!\nabla \!f(\bm{x}^*)]\le0$. Then, given any scalar $0<\gamma<1$, there holds
\begin{align}\label{eq_Vdotexp2}
	\dot V\!\le&- \gamma(\Theta_1+\Theta_2)^{\top}[\nabla \!f(\Theta_1+\Theta_2+\!\bm{x}^*)\!-\!\nabla \!f(\bm{x}^*)]\nonumber\\
	&-\!X^{\top}\! Q^{\top}\!\bar{P}\bar{L}\bar{P}QX\nonumber\\
	=&- \gamma\Theta_1^{\top}[\nabla \!f(\Theta_1+\!\bm{x}^*)\!-\!\nabla \!f(\bm{x}^*)]\nonumber\\
	&- \gamma\Theta_1^{\top}[\nabla \!f(\Theta_1+\Theta_2+\!\bm{x}^*)\!-\!\nabla \!f(\Theta_1+\bm{x}^*)]\nonumber\\
	&- \gamma\Theta_2^{\top}[\nabla \!f(\Theta_1+\Theta_2+\!\bm{x}^*)\!-\!\nabla \!f(\Theta_2+\bm{x}^*)]\nonumber\\
	&- \gamma\Theta_2^{\top}[\nabla \!f(\Theta_2+\!\bm{x}^*)\!-\!\nabla \!f(\bm{x}^*)]\nonumber\\
	&-\!X^{\top}\! Q^{\top}\!\bar{P}\bar{L}\bar{P}QX
\end{align}
Now, from the definition \eqref{eq_defTheta} and the upper-left $0$ entry of equation \eqref{eq_Vdot}, it can be observed that $\bar{L}\bar{P}QX=\bar{L}\bar{P}Q\begin{bmatrix}
	0_{\bar{n}_1}\\X_2
\end{bmatrix}=\bar{L} \Theta_2$. Consequently, one has $\bar{L} \Theta_1=\bar{L}\bar{P}QX-\bar{L}\Theta_2=0$.
Because $\ker \bar{L}=\image \left(\bm{1}_m\otimes I_n\right)$, one has $\Theta_1=\bm{1}_m\otimes q_2$, $ q_2\in\mathbb{R}^{n}$. Recall also that $F(x)=\sum_{i=1}^m f_i(x)$ is strongly convex at $x^*$, thus, there exists a positive $\omega$ such that for any $q_2\in\mathbb{R}^n$,
\begin{align}\label{eq_SCieq}
	&\sum_{i=1}^m\left[q_2^{\top}(\nabla f_i(q_2+x^*)-\nabla f_i(x^*))\right]\nonumber\\
	=~&\left(\bm{1}_m\otimes q_2\right)^{\top}\left(\nabla F(q_2+x^*)-\nabla F(x^*)\right)\geq \omega\|q_2\|_2^2
\end{align}
Thus, for $\Theta_1=\bm{1}_m\otimes q_2$, $q_2\neq0$, 
\begin{align}\label{eq:defconvex}
	& \Theta_1^{\top}[\nabla \!f(\Theta_1+\!\bm{x}^*)\!-\!\nabla \!f(\bm{x}^*)]
	\!=\!\!\sum_{i=1}^m [{q_2}\!^{\top}(\nabla \!f_i(\!q_2\!+\!{x}^*\!)\!-\!\nabla \!f_i(\!{x}^*\!))]\nonumber\\
	&\qquad\ge \omega\|q_2\|_2^2=\frac{\omega}{m}\|\Theta_1\|_2^2
\end{align}
In addition, from \eqref{eq_Vdot} and the fact that $R_2^{\top}\bar{P}\bar{L}\bar{P}R_2$ is positive definite, there exists a positive $\beta_0$ such that 
$$\!X^{\top}\! Q^{\top}\!\bar{P}\bar{L}\bar{P}QX=\!X_2^{\top}\! R_2^{\top}\bar{P}\bar{L}\bar{P}R_2X_2\ge \beta_0\|X_2\|_2^2. $$
From definition \eqref{eq_defTheta} and the fact that the 2-matrix norms of $\bar{P},Q$ are no greater than 1, there holds:
\begin{align}\label{eq_XtoThe}
\!X^{\top}\! Q^{\top}\!\bar{P}\bar{L}\bar{P}QX\ge \beta_0\|\Theta_2\|_2^2
\end{align}
To continue, recall that all $f_i(\cdot)$ are convex, then, they are also locally Lipschitz continuous \cite{D97Nonlinear}. From the convexity, one has 
\begin{align}\label{eq_thetad}
	\Theta_2^{\top}[\nabla \!f(\Theta_2+\!\bm{x}^*)\!-\!\nabla \!f(\bm{x}^*)]\ge0
\end{align}
From the locally Lipschitz continuity and the fact that $X$ is bounded (due to its asymptotic convergence), then there must exist a certain constant $\ell_p>0$ such that
\begin{align}\label{eq_lip}
	\left\|\nabla\! f(\xi_1)\!-\!\nabla\! f(\xi_2)\right\|\le \ell_p\left\|\xi_1-\xi_2\right\|
\end{align}
From \eqref{eq_lip}, by setting $\xi_1=\Theta_1+\Theta_2+\bm{x}^*$, $\xi_2=\Theta_1+\bm{x}^*$; and $\xi_1=\Theta_1+\Theta_2+\bm{x}^*$, $\xi_2=\Theta_2+\bm{x}^*$, one has,
\begin{align}
	&- \!\gamma\Theta_1^{\top}[\nabla \!f(\Theta_1+\Theta_2+\!\bm{x}^*)\!-\!\nabla \!f(\Theta_1+\bm{x}^*)]\le \gamma\ell_p\|\Theta_1\|\|\Theta_2\|\nonumber\\
	&-\! \gamma\Theta_2^{\top}[\nabla \!f(\Theta_1+\Theta_2+\!\bm{x}^*)\!-\!\nabla \!f(\Theta_2+\bm{x}^*)]\le \gamma\ell_p\|\Theta_1\|\|\Theta_2\|\nonumber
\end{align}
Bringing the above equations and \eqref{eq:defconvex}-\eqref{eq_thetad}  into \eqref{eq_Vdotexp2} yields:
\begin{align}\label{eq_Vdotexp4}
	\dot V\!\le&- \frac{\gamma\omega}{m}\|\Theta_1\|_2^2+ 2\gamma\ell_p\|\Theta_1\|\|\Theta_2\|
	-\beta_0\|\Theta_2\|_2^2\nonumber\\
	=&- \left(\sqrt{\frac{\gamma\omega}{2m}}\|\Theta_1\|_2-\sqrt{\frac{2\gamma m \ell_p^{2}}{\omega}}\|\Theta_2\|_2\right)^2\nonumber\\
	&- \frac{\gamma\omega}{2m}\|\Theta_1\|_2^2-(\beta_0-\frac{2\gamma m \ell_p^{2}}{\omega})\|\Theta_2\|_2^2\nonumber\\
	\le&- \frac{\gamma\omega}{2m}\|\Theta_1\|_2^2-\beta_1\|\Theta_2\|_2^2
\end{align}
where $\beta_1=\beta_0-\frac{2\gamma m \ell_p^{2}}{\omega}$. Note that $\beta_0, \ell_p, \omega$ are positive conepants and as in \eqref{eq_Vdotexp2}, $\gamma$ can  be selected arbitrarily from $(0,1)$; hence by choosing sufficiently small $\gamma$, the constant $\beta_1$ can be made strictly positive.
Consider the triangle inequality: 
$$\|\Theta_1\|_2^2+\|\Theta_2\|_2^2\ge \frac{(\|\Theta_1\|_2+\|\Theta_2\|_2)^2}{2}\ge \frac{\|\Theta_1+\Theta_2\|_2^2}{2}$$
and recall the fact that $\Theta_1+\Theta_2=QX$, by choosing $\beta_2=\min (\frac{\gamma\omega}{2m},\beta_1)$, there holds:
\begin{align}\label{eq_Vdotexp5}
	\dot V\!\le&- \frac{\beta_2\|\Theta_1+\Theta_2\|^2}{2}=- \frac{\beta_2\|QX\|_2^2}{2}=- \frac{\beta_2\|X\|_2^2}{2}
\end{align}

To continue, define the following function 
$
V_2(X,Y)=\frac{1}{2}(X+Y)^{\top}(X+Y).
$
Then
\begin{align}\label{eq_dv2}
&\dot{V}_2=(X+Y)^{\top}(\dot{X}+\dot{Y})\nonumber\\
=&-(X+Y)^{\top}[Q^{\top}\bar P[\nabla f(QX\!+\!\bm{x}^*)-\nabla f(\bm{x}^*)]+ Y]\nonumber\\
=&- X^{\top}\!Q^{\top}\![\nabla\! f(QX\!+\!\bm{x}^*)\!-\!\nabla\! f(\bm{x}^*)]
\!-X^{\top}Y\nonumber\\
&- [\nabla\! f(QX\!+\!\bm{x}^*)\!-\!\nabla\! f(\bm{x}^*)]^{\top}Y -Y^{\top}Y \nonumber\\
=&- X^{\top}\!Q^{\top}\![\nabla\! f(QX\!+\!\bm{x}^*)\!-\!\nabla\! f(\bm{x}^*)]-\frac{1}{2}Y^{\top}Y\nonumber\\
&-\left\|X+\frac{1}{2}Y\right\|_2^2-\left\|[\nabla\! f(QX\!+\!\bm{x}^*)\!-\!\nabla\! f(\bm{x}^*)]+\frac{1}{2}Y\right\|_2^2\nonumber\\
&+\left\|X\right\|_2^2+\left\|\nabla\! f(QX\!+\!\bm{x}^*)\!-\!\nabla\! f(\bm{x}^*)\right\|_2^2\nonumber\\
\le &-\frac{1}{2}Y^{\top}Y+ \left\|X\right\|_2^2+\left\|\nabla\! f(QX\!+\!\bm{x}^*)\!-\!\nabla\! f(\bm{x}^*)\right\|_2^2\nonumber\\
\le &-\frac{1}{2}\left\|Y\right\|_2^2+ (1+\ell_p^2)\left\|X\right\|_2^2
\end{align}
The last inequality holds because of the inequality associated with Lipschitz continuity given in \eqref{eq_lip}. Now consider 
\begin{align}\label{eq_defwtV}
	\widetilde{V}=&\frac{2(2+\ell_p^2)}{\beta_2}V+V_2.
\end{align}
From equations \eqref{eq_Vdotexp5} and \eqref{eq_dv2}, there holds:
\begin{align}\label{eq_dV2}
	\dot{\widetilde{V}}\le&-\left\|X\right\|_2^2-\frac{1}{2} \left\|Y\right\|_2^2
\end{align}
By the definitions of $\widetilde{V}$, ${V}$, and ${V_2}$, obviously, there exists a positive constant $\beta$ such that $\left\|X\right\|_2^2+\frac{1}{2} \left\|Y\right\|_2^2\ge\beta \widetilde{V}$. Therefore,
\begin{align}\label{eq_dV2f}
	\dot{\widetilde{V}}\le-\beta \widetilde{V}
\end{align}
This established the global exponential convergence of the system and completes the proof the theorem. \qed

\vspace{-0.5em}
\subsection{Proof of Corollary \ref{Cl1}}
Consider the proof of Theorem \ref{thm1} with replacement of  equation \eqref{Dupdatex}  by equation \eqref{DupdatexD}, and define $\bm{v}=\col\{v_1,\cdots,v_m\}$. From the definition of $\widetilde{V}$, and  the obtained $\dot{\widetilde{V}}$ in equation \eqref{eq_dV2f}, there must exist positive constants  $\eta_1$, $\eta_2$, $\eta$, such that

\begin{align}
	\dot{\widetilde{V}}&\le-\beta \widetilde{V}+\eta_1\|\bm{v}\|_2\|X\|_2+\eta_2\|\bm{v}\|_2\|Y\|_2 \nonumber\\
	&\le-\beta \widetilde{V}+\eta\|\bm{v}\|_2\sqrt{\widetilde{V}}
\end{align}
Further since $\|\bm{v}\|_2$ is bounded, i.e. $\|\bm{v}\|_2\le \epsilon$ for some positive $\epsilon$, then as $t\to \infty$, one has 
\begin{align}
	\widetilde{V}\le \frac{\eta^2\epsilon^2}{\beta^2}.
\end{align}
is bounded. It follows that $\|x_i(t)-x^*\|$ is bounded for all $i=1,\cdots,m$. This completes the proof. \qed

\section{Simulation}
\subsection{The exponential convergence rate}\label{sec_Expo}
First consider an example with $m=5$, $n=20$, and $n_i=3$, for $i=1,\cdots,5$. Suppose the agents in the network have the following neighbor relations: $\mathcal{N}_1=\{1,2,3,4\}$, $\mathcal{N}_2=\{1,2,3\}$, $\mathcal{N}_3=\{1,2,3,4\}$, $\mathcal{N}_4=\{1,3,4,5\}$, $\mathcal{N}_5=\{4,5\}$, which ensures the network is undirected and connected. 
We let each agent $i$ know a local objective function $f_i(x)$ and a local constraint $A_ix=b_i$ such that
\begin{align}\label{eq_Simu_Setups}
&f_1(x)=\|x\|_2^2 \qquad\qquad f_2(x)=\|x-c_2\|_2^2  \nonumber\\
&f_3(x)=\sum_{k=1}^{20} e^{x[k]}  \qquad\quad f_4(x)=\sum_{k=1}^{20} e^{-2x[k]} \\
&f_5(x)=\|x-c_5\|_2^4  \nonumber
\end{align}
where $x[k]$ denotes the $k$th entry of vector $x$; $A_i\in\mathbb{R}^{3\times 20}$, $b_i\in\image A_i\subset \mathbb{R}^{3}$ and $c_2,c_5\in\mathbb{R}^{20}$ are constant matrices/vectors. 
We let the simulation configuration satisfy Assumption \ref{ass2}.

\begin{figure}[h]
	\centering
	\includegraphics[width=8.2 cm]{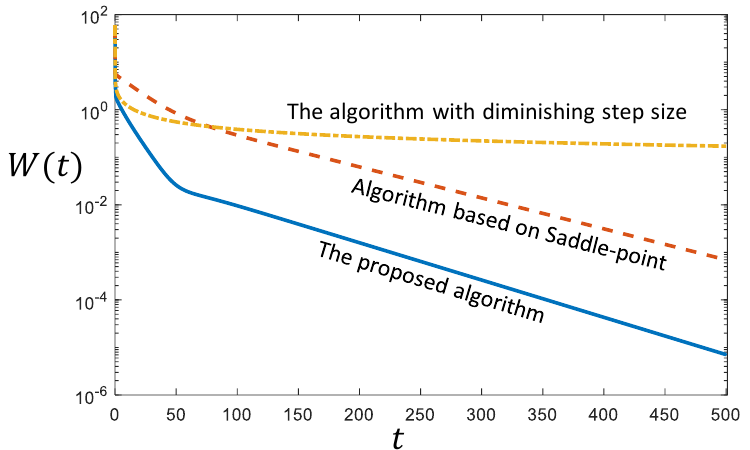}
	\caption{The 5-agent case, exponential convergence of the algorithm, with simulation execution time 0.13s.}
	\label{Fig_expo}
\end{figure}

In order to validate Theorem \ref{thm1}, we let each agent initialize its local state $x_i\in\mathbb{R}^{20}$ as $A_ix_i(0)=b_i$ and then update its state by equation \eqref{Dupdatex}.  Define the following function:
\begin{align}
W(t)=\sum_{i=1}^5 \|x_i(t)-x^*\|_2^2,
\end{align}
for which $W(t)=0$ if and only if all $x_i(t)=x^*$ for all $i=1,\cdots,5$,
where $x^*$  is the unique minimizer of $F(x)=\sum_{i=1}^m f_i(x)$ subject to $A_ix=b_i$, $i=1\cdots,5$. The simulation result is obtained using the Ode45 solver of MATLAB, with a computer equipped with Intel 6700 CPU. The result is presented in Fig. \ref{Fig_expo}, where the Y-axis is scaled by $\log(\cdot)$. The constant slope of the curve $W(t)$ (with $W(t)$ converging to $0$) indicates the exponential convergence of the algorithm, which validates Theorem \ref{thm1}. In contrast, algorithms with diminishing gain ($\alpha(t)=1/t$) only achieve an asymptotic convergence rate. Note that the continuous algorithms based on  saddle-point dynamics (primal-dual)\cite{BJ14TAC} can also achieve an exponential convergence rate, but as discussed earlier, require extra states that are transferred across the network. 
\noindent

\vspace{-1em}
\subsection{Robustness towards disturbance}
Here, we investigate the impact of bounded disturbance  to the proposed algorithm. Consider a connected, undirected network with $m=30$ agents. We let each agent $i$ know a local state $x_i\in\mathbb{R}^5$, a local objective function $f_i(x_i)$ and a local constraint $A_ix=b_i$. 
Note that the local objective functions basically follow one of the linear, norm, or exponential forms, which are similar to the examples provided in \eqref{eq_Simu_Setups}. The coefficients of the objective functions are randomly chosen, and we make sure Assumptions \ref{ass2} holds.
In addition, for each agent, we use MATLAB Random Source tool to introduce a bounded disturbance $v_i$, such that each entry of $v_i$ is chosen from $[0,~0.01]$ with sampling time 0.1s. 
We compare the performance of update  \eqref{DupdatexD} and the update
\begin{align}\label{eq_Alg_optdis}
	\dot{x}_i= -P_i\left(\alpha(t)\nabla f_i(x_i)+\sum_{j\in\mathcal{N}_i} (x_i -x_j)\right)+v_i,
\end{align}
equipped with diminishing gain $\alpha(t)=1/t$ (the discrete-time versions of this update are proposed in \cite{AAP10TAC,PWY16Auto}).
\noindent\begin{figure}[h]
	\centering
	\includegraphics[width=8.2 cm]{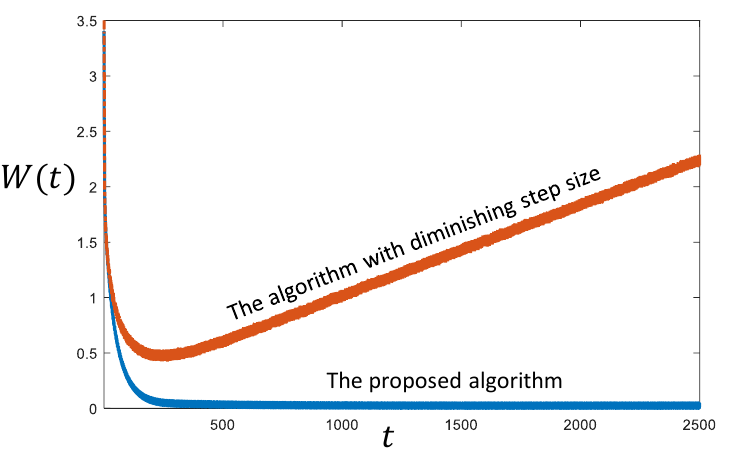}
	\caption{The 30-agent case, under bounded disturbance, with simulation execution time 1.05s.}
	\label{Fig_noise}
\end{figure}

For both algorithms, we use the same initial state and the Ode45 solver of MATLAB to perform simulation. The obtained error curves $W(t)$ are shown in Fig. \ref{Fig_noise}, which is a representative example we choose from the many random simulations we performed.
It can be observed that for the update with diminishing gain, the non-zero disturbance will accumulate with time $t$ and finally lead the curve $W(t)$ to  blow up. This means the agents' states are not able to converge to the optimum point $x^*$.  
For the proposed update, the curve of $W(t)$ does not grow with time $t$, instead, it converges to a bounded neighborhood of  $x^*$, which validates Corollary \ref{Cl1}.

\section{Conclusion}
In this paper, by incorporating the idea of integral feedback, we proposed a continuous-time distributed algorithm which is able to solve a constrained distributed optimization problem with exponential convergence rate. 
To sum up, the proposed algorithm (a) does represent an increase in the state dimension at each agent over an algorithm with diminishing gain (b) in comparison with other algorithms which use an increase of dimension to avoid diminishing gain, does not impose an additional burden on the communication bandwidth (c) presents relaxed condition on objective function to guarantee exponentially fast convergence, (d) offers robustness against disturbance.
Future work includes the generalization of the proposed algorithm to discrete-time update; to
time-varying directed networks;  and application to general local constraints other than those expressed by a local linear equation.

\bibliographystyle{IEEEtran}

\bibliography{Ref,Shaoshuai,Linear_nonlinear_equations}  

\end{document}